\documentclass[12pt]{article}
\usepackage{amsmath,amsfonts,amssymb,amscd}

\newtheorem{thm}{Theorem}[section]

\newtheorem{lem}[thm]{Lemma}

\newtheorem{rem}[thm]{Remark}

\makeatletter \@addtoreset{figure}{section} \makeatother
\makeatletter
\long\def\@makecaption#1#2{%
   \vskip 10\p@
   \setbox\@tempboxa\hbox{{#1}\ \ #2}%
   \ifdim \wd\@tempboxa >\hsize
       {#1}\ \ #2\par
   \else
       \hbox to\hsize{\hfil\box\@tempboxa\hfil}%
   \fi}
\makeatother

\def\qed{\hfill \rule{4pt}{7pt}}
\def\pf{\noindent {\it Proof.} }
\title{\bf Long heterochromatic paths in heterochromatic triangle free graphs
\footnote{Supported by NSFC, PCSIRT and the ``973" program.}}
\author{\normalsize\itshape He Chen and Xueliang Li\\
\normalsize\itshape Center for Combinatorics and LPMC-TJKLC,\\
\normalsize\itshape Nankai University, Tianjin 300071, P.R. China\\
\normalsize chenhe@mail.nankai.edu.cn, lxl@nankai.edu.cn}

\date{}

\begin{document}

\maketitle

\begin{abstract}
In this paper, graphs under consideration are always edge-colored.
We consider long heterochromatic paths in heterochromatic triangle
free graphs. Two kinds of such graphs are considered, one is
complete graphs with Gallai colorings, i.e., heterochromatic
triangle free complete graphs; the other is heterochromatic triangle
free graphs with $k$-good colorings, i.e., minimum color degree at
least $k$. For the heterochromatic triangle free graphs $K_n$, we
obtain that for every vertex $v\in V(K_n)$, $K_n$ has a
heterochromatic $v$-path of length at least $d^c(v)$; whereas for
the heterochromatic triangle free graphs $G$ we show that if, for
any vertex $v\in V(G)$, $d^c(v)\geq k\geq 6$, then $G$ a
heterochromatic
path of length at least $\frac{3k}{4}$. \\[3mm]
\\[3mm]

\end{abstract}

\section{Introduction}

We use Bondy and Murty \cite{B-M} for terminology and notations not
defined here and consider simple graphs only.

Let $G=(V,E)$ be a graph. By an {\it edge coloring} of $G$ we will
mean a function $C: E\rightarrow \mathbb{N} $, the set of natural
numbers. If $G$ is assigned such a coloring, then we say that $G$ is
an {\it edge-colored graph}. Denote the edge-colored graph by
$(G,C)$, and call $C(e)$ the {\it color} of the edge $e\in E$. We
say that $C(uv)=\emptyset$ if $uv\notin E(G)$ for $u,v\in V(G)$. For
a subgraph $H$ of $G$, we denote $C(H)=\{C(e) \ | \ e\in E(H)\}$ and
$c(H)=|C(H)|$. For a vertex $v$ of $G$, we say that color $i$ is
{\it presented at} vertex $v$ if some edge incident with $v$ has
color $i$. The {\it color degree $d^c(v)$} is the number of
different colors that are presented at $v$, and the {\it color
neighborhood} $CN(v)$ is the set of different colors that are
presented at $v$. All graphs considered in this paper are
edge-colored. For a positive integer $k$, a coloring of a graph is
called {\it $k$-good} if the minimum color degree of the graph is at
least $k$. A path, or a cycle, or any subgraph is called {\it
heterochromatic (rainbow, or multicolored)} if any two edges of it
have different colors. A graph is called {\it heterochromatic
triangle free} if it does not contain any (induced) heterochromatic
triangles. If $u$ and $v$ are two vertices on a path $P$, $uPv$ will
denote the segment of $P$ from $u$ to $v$, whereas $vP^{-1}u$ will
denote the same segment but from $v$ to $u$. A path is called a {\it
$v$-path} if it starts from the vertex $v$.

There are a lot of existing literature dealing with the existence of
paths and cycles with special properties in edge-colored graphs. The
heterochromatic Hamiltonian cycle or path problem was studied by
Hahn and Thomassen \cite{H-T}, R\"{o}dl and Winkler (see
\cite{F-R}), Frieze and Reed \cite{F-R}, and Albert, Frieze and Reed
\cite{A-F-R}. In \cite{A-J-Z}, Axenovich, Jiang and Tuza gave the
range of the maximum $k$ such that there exists a $k$-good coloring
of $E(K_n)$ that contains no properly colored copy of a path with
fixed number of edges, no heterochromatic copy of a path with fixed
number of edges, no properly colored copy of a cycle with fixed
number of edges and no heterochromatic copy of a cycle with fixed
number of edges, respectively. In \cite{E-T-1}, Erd\"{o}s and Tuza
studied the heterochromatic paths in infinite complete graph
$K_\omega$. In \cite{E-T-2}, Erd\"{o}s and Tuza studied the values
of $k$, such that every $k$-good coloring of $K_n$ contains a
heterochromatic copy of $F$ where $F$ is a given graph with $m$
edges ($m<n/k$). In \cite{M-S-T}, Manoussakis, Spyratos and Tuza
studied $(s,t)$-cycles in $2$-edge-colored graphs, where an
$(s,t)$-cycle is a cycle of length $s+t$ and $s$ consecutive edges
are in one color and the remaining $t$ edges are in the other color.
In \cite{M-S-T-V}, Manoussakis, Spyratos, Tuza and Voigt studied
conditions on the minimum number $k$ of colors, sufficient for the
existence of given types (such as families of internally pairwise
vertex-disjoint paths with common endpoints, Hamiltonian paths and
Hamiltonian cycles, cycles with a given lower bound of their length,
spanning trees, stars, and cliques ) of properly edge-colored
subgraphs in a $k$-edge-colored complete graph. In \cite{C-M-M},
Chou, Manoussakis, Megalaki, Spyratos and Tuza showed that for a
2-edge-colored graph $G$ and three specified vertices $x, y$ and
$z$, to decide whether there exists a color-alternating path from
$x$ to $y$ passing through $z$ is $NP$-complete. Many results in
these mentioned papers were proved by using probabilistic methods.

In \cite{A-J-Z}, Axenovich, Jiang and Tuza considered the local
variation of anti-Ramsey problem. Namely, they studied the maximum
integer $k$, denoted by $g(n,H)$, such that there exists a $k$-good
edge coloring of $K_n$ that does not contain any heterochromatic
copy of a given graph $H$. They showed that for a fixed integer
$k\geq 2$, $k-1\leq g(n, P_{k+1})\leq 2k-3$, i.e., if $K_n$ is
edge-colored by a $(2k-2)$-good coloring, then there must exist a
heterochromatic path $P_{k+1}$, there exists a $(k-1)$-good coloring
of $K_n$ such that no heterochromatic path $P_{k+1}$ exists.

In \cite{B-L}, the authors considered the long heterochromatic paths
in general graphs with a $k$-good coloring and showed that if $G$ is
an edge-colored graph with $d^c(v)\geq k$ (color degree condition)
for every vertex $v$ of $G$, then $G$ has a heterochromatic $v$-path
of length at least $\lceil\frac{k+1}{2}\rceil$. In
\cite{C-L-1,C-L-2}, we got some better bound of the length of
longest heterochromatic paths in general graphs with a $k$-good
coloring.

\begin{thm}\label{k leq 7}\cite{C-L-1}
Let $G$ be an edge-colored graph and $3\leq k\leq 7$ an integer.
Suppose that $d^c(v)\geq k$ for every vertex $v$ of $G$. Then $G$
has a heterochromatic path of length at least $k-1$.
\end{thm}

\begin{thm}\label{2/3 k}\cite{C-L-2}
Let $G$ be an edge-colored graph. If $d^c(v)\geq k\geq 7$ for any
vertex $v\in V(G)$, then $G$ has a heterochromatic path of length at
least $\lceil\frac{2k}{3}\rceil+1$.
\end{thm}

In \cite{C-L-3}, we showed that if $|CN(u)\cup CN(v)|\geq s$ (color
neighborhood union condition) for every pair of vertices $u$ and $v$
of $G$, then $G$ has a heterochromatic path of length at least
$\lceil\frac{s+1}{2}\rceil$, and gave examples to show that the
lower bound is best possible in some sense.

Some special edge colorings have also been studied, such as {\it
Gallai colorings}, which is defined to be the edge colorings of
complete graphs in which no heterochromatic triangles exist. In
\cite{G-S}, Gy\'{a}rf\'{a}s and Simonyi studied the existence of
special monochromatic spanning trees in such colorings, they also
determined the size of largest monochromatic stars guaranteed to
occur.

In this paper, we consider long heterochromatic paths in complete
graphs $K_n$ with Gallai colorings, i.e., heterochromatic triangle
free complete graphs, and in heterochromatic triangle free graphs
with $k$-good colorings. We obtain that if $K_n$ is heterochromatic
triangle free, then for every vertex $v\in V(K_n)$, $K_n$ has a
heterochromatic $v$-path of length at least $d^c(v)$. For the
heterochromatic triangle free general graphs $G$, we show that if
$d^c(v)\geq k\geq 6$ for any vertex $v\in V(G)$, then $G$ has a
heterochromatic path of length at least $\frac{3k}{4}$.

\section{Heterochromatic triangle free complete graphs}

In this section, we consider a heterochromatic triangle free
complete graph $G$, and try to find a long heterochromatic path from
it.

\begin{thm}\label{triangle-free-gallai-thm-1}
Suppose $G$ is a heterochromatic triangle free complete graph. Then
for every vertex $u$ in $G$, $G$ has a heterochromatic $u$-path of
length at least $d^c(u)$.
\end{thm}
\pf Let $u$ be any vertex of $G$ and let $d^c(u)=k$. Suppose
$v_1,v_2,\ldots, v_k$ are $k$ different neighbors of $u$ such that
the $k$ edges $uv_1,uv_2,\ldots, uv_k$ all have distinct colors (see
Figure \ref{fig-triangle-free-1}).
\begin{figure}[h,t]
\begin{center}
\begin{picture}(100,55)
\put(10,50){\circle{2}} \put(30,50){\circle{2}}
\put(70,50){\circle{2}} \put(90,50){\circle{2}}
\put(50,10){\circle{2}} \multiput(40,50)(4,0){6}{\circle*{1}}
\put(50,10){\line(-1,1){40}} \put(50,10){\line(-1,2){20}}
\put(50,10){\line(1,2){20}} \put(50,10){\line(1,1){40}}
\put(6,53){$v_1$} \put(26,53){$v_2$} \put(66,53){$v_{k-1}$}
\put(92,53){$v_k$} \put(47,2){$u$}
\end{picture}
\caption{} \label{fig-triangle-free-1}
\end{center}
\end{figure}

At first, we will construct a path $P$ by the following inductive
algorithm.

{\bf Algorithm}

{\bf Step 1.} If the two edges $v_1v_2$ and $uv_1$ have the same
color, we let $w_1=v_1$, $w_2=v_2$; otherwise, we let $w_1=v_2$,
$w_2=v_1$. Set $P_1=w_1w_2$.

If {\bf Step $i-1$} is finished and we have obtained the path
$P_{i-1}=w_1w_2\ldots w_i$. Then

{\bf Step $i$.} If the two edges $w_iv_{i+1}$ and $uw_i$ have the
same color, we let $w_{i+1}=v_{i+1}$.

Otherwise, if the color of the edge $w_{i-1}v_{i+1}$ is the same as
the color of the edge $uw_{i-1}$, we let $w_{i+1}=w_i$ and
$w_i=v_{i+1}$.

Otherwise, let $j_0$ be the maximum integer $j$ such that the two
edges $w_jv_{i+1}$ and $uv_{i+1}$ have the same color and the colors
of the two edges $w_{j-1}v_{i+1}$ and $uv_{i+1}$ are distinct. If
all the $i$ edges
$w_1v_{i+1},w_2v_{i+1},\ldots,w_{i-1}v_{i+1},uv_{i+1}$ have the same
color, we set $j_0=1$.

If $j_0=1$, let $w_{i+1}=w_i$, $w_i=w_{i-1}$, $w_{i-1}=w_{i-2}$,
$\ldots$, $w_3=w_2$, $w_2=w_1$, $w_1=v_{i+1}$. Otherwise, $2\leq
j_0\leq i-1$, let $w_{i+1}=w_i$, $w_i=w_{i-1}$, $\ldots$,
$w_{j_0+1}=w_{j_0}$, $w_{j_0}=v_{i+1}$. Set $P_i=w_1w_2\ldots
w_{i+1}$.

Continue the process till $i=k$ and we obtain the path $P=P_k$.

Then, we will prove the following claim about the path $P$ obtained
from the algorithm.

\noindent {\bf Claim.} The path $P=w_1w_2\ldots w_k$ obtained from
the algorithm is heterochromatic, the vertex set of path $P$ is
actually the set $\{v_1,v_2,\ldots,v_k\}$, and for each $1\leq l\leq
k-1$, the two edges $w_lw_{l+1}$ and $uw_l$ have the same color.

\pf To prove the claim, we will show that for each $i$ ($1\leq i\leq
k-1$), the path $P_i=w_1w_2\ldots w_{i+1}$ we obtained after Step
$i$ satisfies that the vertex set of $P_i$ is actually the set
$\{v_1,v_2,\ldots,v_{i+1}\}$, and for each $1\leq l\leq i$, the two
edges $w_lw_{l+1}$ and $uw_l$ have the same color.

When $i=1$, since there is no heterochromatic triangle in $K_n$, the
edge $v_1v_2$ has the same color as the color of the edge $uv_1$ or
$uv_2$. From Step 1, we can easily see that no matter which color
the edge $v_1v_2$ has, the path $P=w_1w_2$ that is obtained after
Step 1 contains actually two vertices $v_1$ and $v_2$, and the the
two edges $w_1w_2$ and $uw_1$ have the same color.

Suppose Step $i-1$ has been finished, and the path $P=w_1w_2\ldots
w_i$ we have obtained now contains actually $i$ vertices
$v_1,v_2,\ldots,v_i$, and for each $1\leq l\leq i-1$, the two edges
$w_lw_{l+1}$ and $uw_l$ have the same color, which is named to be
color $c_l$. See Figure \ref{fig-triangle-free-2}. Now we consider
the path $P_i=w_1w_2\ldots w_{i+1}$ we obtained after Step $i$.

\begin{figure}[h,t]
\begin{center}
\begin{picture}(300,100)
\put(10,90){\circle{2}} \put(50,90){\circle{2}}
\put(90,90){\circle{2}} \put(170,90){\circle{2}}
\put(210,90){\circle{2}} \put(290,90){\circle{2}}
\put(130,10){\circle{2}} \multiput(100,90)(4,0){16}{\circle*{1}}
\put(130,10){\line(-3,2){120}} \put(130,10){\line(-1,1){80}}
\put(130,10){\line(-1,2){40}} \put(130,10){\line(1,2){40}}
\put(130,10){\line(1,1){80}} \put(130,10){\line(2,1){160}}
\put(10,90){\line(1,0){40}} \put(50,90){\line(1,0){40}}
\put(170,90){\line(1,0){40}} \put(210,90){\line(1,0){80}}
\put(6,93){$w_1$} \put(46,93){$w_2$} \put(86,93){$w_3$}
\put(162,93){$w_{i-1}$} \put(212,93){$w_i$} \put(283,93){$v_{i+1}$}
\put(127,2){$u$} \put(27,88){$c_1$} \put(62,55){$c_1$}
\put(67,88){$c_2$} \put(87,55){$c_2$} \put(107,55){$c_3$}
\put(187,88){$c_{i-1}$} \put(140,55){$c_{i-1}$} \put(180,55){$c_i$}
\put(227,55){$c_{i+1}$}
\end{picture}
\caption{} \label{fig-triangle-free-2}
\end{center}
\end{figure}

For convenience, if the path $P_i=w_1w_2\ldots w_{i+1}$ satisfies
that the vertex set of $P_i$ is actually the set
$\{v_1,v_2,\ldots,v_{i+1}\}$, and for each $1\leq l\leq i$, the two
edges $w_lw_{l+1}$ and $uw_l$ have the same color, we say that $P_i$
satisfies {\bf Condition A}.

If the two edges $w_iv_{i+1}$ and $uw_i$ have the same color, then
$w_{i+1}=v_{i+1}$. So, $P_i$ obviously satisfies Condition A.

Otherwise, since the triangle $uw_iwv_{i+1}$ is not heterochromatic,
the two edges $w_iv_{i+1}$ and $uv_{i+1}$ have the same color. In
this case, if the color of the edge $w_{i-1}v_{i+1}$ is the same as
the color of the edge $uw_{i-1}$, then $w_{i+1}=w_i$ and
$w_i=v_{i+1}$ by Step $i$. Thus, $P_i$ satisfies Condition A.

Now we consider the case when the two edges $w_iv_{i+1}$ and
$uv_{i+1}$ have the same color, and the two edges $w_{i-1}v_{i+1}$
and $uw_{i-1}$ have two distinct colors. Noticing that the triangle
$uw_{i-1}v_{i+1}$ is not heterochromatic, we can conclude that the
two edges $w_{i-1}v_{i+1}$ and $uv_{i+1}$ have the same color.

If all the $i$ edges
$w_1v_{i+1},w_2v_{i+1},\ldots,w_{i-1}v_{i+1},uv_{i+1}$ have the same
color, we have that the two edges $w_1v_{i+1}$ and $uv_{i+1}$ have
the same color and for each $l$ ($1\leq l\leq i-1$), the two edges
$uw_l$ and $w_lw_{l+1}$ have the same color. Thus, the path
$v_{i+1}w_1w_2\ldots w_i$ satisfies Condition A. This implies that
when we set $w_{i+1}=w_i$, $w_i=w_{i-1}$, $w_{i-1}=w_{i-2}$,
$\ldots$, $w_3=w_2$, $w_2=w_1$, $w_1=v_{i+1}$, the path
$P_i=w_1w_2\ldots w_iw_{i+1}$ satisfies Condition A.

Otherwise, we can find a maximum integer $j$, say $j_0$, such that
the two edges $w_jv_{i+1}$ and $uv_{i+1}$ have the same color, and
the colors of the two edges $w_{j-1}v_{i+1}$ and $uv_{i+1}$ are
distinct. It is clear from the discussion above that in this case,
$2\leq j_0\leq i-1$. Then the vertex $w_{j_0-1}\in
\{v_1,v_2,\ldots,v_{i}\}$. By the assumption that the colors that
the $k$ edges $uv_1,uv_2,\ldots,uv_k$ all have distinct edges, we
have that the two edges $uw_{j_0-1}$ and $uv_{i+1}$ have two
distinct colors. On the other hand, the two edges $w_{j_0-1}v_{i+1}$
and $uv_{i+1}$ have two distinct colors. So we can conclude that the
two edges $w_{j_0-1}v_{i+1}$ and $uw_{j_0-1}$ have the same color,
because the triangle $uw_{j_0-1}v_{i+1}$ is not heterochromatic.
Then we can conclude that the path $w_1w_2\ldots
w_{j_0-1}v_{i+1}w_{j_0}w_{j_0+1}\ldots w_i$ satisfies Condition A.
This implies that if we set $w_{i+1}=w_i$, $w_i=w_{i-1}$, $\ldots$,
$w_{j_0+1}=w_{j_0}$, $w_{j_0}=v_{i+1}$, the path $P_i=w_1w_2\ldots
w_iw_{i+1}$ satisfies Condition A. The Claim is thus proved by
induction. \qed

Now we turn back to the proof of the theorem.

Since the path $P$ we obtained from the algorithm satisfies all the
conditions in the Claim, the color of the edge $w_ku$ does not
appear on the path $P$, and so the path $w_1Pw_ku$ is a
heterochromatic $u$-path of length $k$. The proof is thus complete.
\qed

Using the theorem above, we can easily get the following result as a
corollary.
\begin{thm}\label{triangle-free-gallai-thm-2}
Suppose $G$ is a heterochromatic triangle free complete graph. If
the maximum color degree among all the vertices in $G$ is $k$, i.e.,
$\max_{v\in V(G)}{d^c(v)}=k$, then there is a heterochromatic path
of length at least $k$ in $G$.
\end{thm}

\section{Heterochromatic triangle free general graphs}

In this section, we consider long heterochromatic paths in a
heterochromatic triangle free general graph. Before we give our main
theorem, we would like to give some properties about this special
kind of edge colored graphs.

\begin{lem}\label{triangle-free-lem}
Suppose $G$ is a heterochromatic triangle free graph, and
$P=u_0u_1\ldots u_l$ is a heterochromatic path of length $l\geq 5$.
If the two edges $u_0u_i$ and $u_0u_j$ ($2\leq i<i+1<j\leq l$) exist
and their colors are distinct and do not appear on the path $P$,
then either there exists an integer $s$, $i<s<j$, such that the edge
$u_0u_s$ does not exist, or there exist two integers $s$ and $t$
($i< s<t\leq j$), such that the two edges $u_0u_s$ and $u_0u_t$ have
the same color.
\end{lem}
\pf We will prove it by contradiction.

Suppose that we cannot get the conclusion, which implies that all
the edges $u_0u_{i+1}, u_0u_{i+2},\ldots, u_0u_{j-1}, u_0u_j$ exist
in $G$ and they all have distinct colors.

First, we consider the triangle $u_0u_{j-1}u_j$. Since it is not
heterochromatic, the color of the edge $u_0u_j$ does not appear on
the path $P$, and the two edges $u_0u_{j-1}$ and $u_0u_j$ have two
distinct colors, we have that the two edges $u_0u_{j-1}$ and
$u_{j-1}u_j$ have the same color.

Now we consider the triangle $u_0u_{j-2}u_{j-1}$. As the two edges
$u_0u_{j-1}$ and $u_{j-1}u_j$ have the same color, and the path $P$
is heterochromatic, we have that the two edges $u_0u_{j-1}$ and
$u_{j-2}u_{j-1}$ have two distinct colors. On the other hand, by the
assumption, the triangle $u_0u_{j-2}u_{j-1}$ is not heterochromatic,
and the two edges $u_0u_{j-2}$ and $u_0u_{j-1}$ have two distinct
colors. So the two edges $u_0u_{j-2}$ and $u_{j-2}u_{j-1}$ have the
same color.

In the same way, we can get, orderly, the edge $u_0u_{j-3}$ has the
same color as the edge $u_{j-3}u_{j-2}$ has, $\ldots$, the edge
$u_0u_{i+1}$ has the same color as the edge $u_{i+1}u_{i+2}$ has.
Then the triangle $u_0u_iu_{i+1}$ is heterochromatic, a
contradiction, which completes the proof. \qed

In a similar way, we can get the following property.

\begin{lem}\label{triangle-free-lem-1}
Suppose $G$ is a heterochromatic triangle free graph, and
$P=u_0u_1\ldots u_l$ is a heterochromatic path of length $l\geq 5$.
If the edge $u_0u_i$ exists and the color of it does not appear on
the path $P$, then $i\geq 3$, and either there exists an integer
$s$, $2\leq s<i$, such that the edge $u_0u_s$ does not exist, or
there exist two integers $s$ and $t$ ($2\leq s<t\leq i$), such that
the two edges $u_0u_s$ and $u_0u_t$ have the same color.
\end{lem}

Now we can state our main theorem.
\begin{thm}
Suppose $G$ is a heterochromatic triangle free graph. If $d^c(v)\geq
k\geq 6$ for any vertex $v\in V(G)$, then $G$ has a heterochromatic
path of length at least $\frac{3k}{4}$.
\end{thm}
\pf Suppose $P=u_0u_1u_2\ldots u_l$ is one of the longest
heterochromatic paths in $G$. Assume that $CN(u_0)$ has $s$
different colors not appearing on $P$, and $CN(u_l)$ has $t$
different colors not appearing on $P$. Then there exist $s$
different vertices $u_{x_1},u_{x_2},\ldots,u_{x_s}$ on the path $P$,
where $2\leq x_1<x_2<\ldots <x_s\leq l$, such that the colors of the
$s$ edges $u_0u_{x_1},u_0u_{x_2},\ldots,u_0u_{x_s}$ are all distinct
and do not appear on $P$. There also exist $t$ different vertices
$u_{y_1},u_{y_2},\ldots,u_{y_t}$ on the path $P$, where $0\leq
y_1<y_2<\ldots y_t\leq l-2$, such that the colors of the $t$ edges
$u_{y_1}u_l,u_{y_2}u_l,\ldots,u_{y_t}u_l$ are all distinct and do
not appear on $P$. Since there exists no heterochromatic triangle in
$G$, we have $x_1\geq 3$, $x_{i+1}> x_i+1$ for $i=1,2,\ldots, s-1$;
$y_t\leq l-3$, $y_{j+1}>y_j+1$ for $j=1,2,\ldots, t-1$.

Since $k\geq 6$, we can conclude from Theorems \ref{k leq 7} and
\ref{2/3 k} that the path $P$ is of length $l\geq 5$. By Lemma
\ref{triangle-free-lem-1}, we have that $$|\{C(u_0u_2),
C(u_0u_3),\ldots,C(u_0u_{x_1})\}|\leq x_1-2.$$ We can also get from
Lemma \ref{triangle-free-lem} that for any $1\leq i\leq s-1$,
$$|\{C(u_0u_{x_i+1}),C(u_0u_{x_i+2}),\ldots, C(u_0u_{x_{i+1}-1}),
C(u_0u_{x_{i+1}})\}|\leq x_{i+1}-x_i-1.$$ So
\begin{eqnarray}\label{ineq-1}
\begin{array}{ll}
 & |\{C(u_0u_1), C(u_0u_2),\ldots,C(u_0u_{l-1}), C(u_0u_l)\}|\\
\leq & |\{C(u_0u_1))\}|+|\{C(u_0u_2), C(u_0u_3),\ldots,C(u_0u_{x_1})\}|\\
 & + |\{C(u_0u_{x_1+1}),C(u_0u_{x_1+2}),\ldots, C(u_0u_{x_2-1}),
C(u_0u_{x_2})\}|\\
 & + |\{C(u_0u_{x_2+1}),C(u_0u_{x_2+2}),\ldots, C(u_0u_{x_3-1}),
C(u_0u_{x_3})\}|\\
 & + \ldots\\
 & + |\{C(u_0u_{x_{s-1}+1}),C(u_0u_{x_{s-1}+2}),\ldots, C(u_0u_{x_s-1}),
C(u_0u_{x_s})\}|\\
 & + |\{C(u_0u_{x_s+1}),\ldots, C(u_0u_{l-1}), C(u_0u_l)\}|\\
\leq & 1+(x_1-2)+(x_2-x_1-1)+(x_3-x_2-1)+\ldots+(x_s-x_{s-1}-1)+(l-x_s)\\
= & l-s.
\end{array}
\end{eqnarray}

On the other hand, for any vertex $v$ which is adjacent to $u_0$ but
does not belong to the path $P$, the color of the edge $u_0v$ is not
same as the color of the edge $u_{y_j}u_{y_j+1}$ for any $1\leq
j\leq t$, for otherwise, $vu_0Pu_{y_j}u_lP^{-1}u_{y_j+1}$ is a
heterochromatic path of length $l+1$, a contradiction. So we have
$CN(u_0)\setminus \{C(u_0u_i): 1\leq i\leq l\}\subseteq
C(P)\setminus \{C(u_{y_j}u_{y_j+1}):1\leq j\leq t\}$, and then
\begin{eqnarray}\label{ineq-2}
|CN(u_0)\setminus \{C(u_0u_i): 1\leq i\leq l\}|\leq l-t.
\end{eqnarray}

From Inequalities \ref{ineq-1} and \ref{ineq-2}, we have
\begin{eqnarray}\label{ineq-3}
\begin{array}{lll}
k & \leq & |CN(u_0)|\\
&\leq & |CN(u_0)\setminus \{C(u_0u_i): 1\leq i\leq
l\}|+|\{C(u_0u_i): 1\leq i\leq l\}|\\
&\leq & (l-t)+(l-s)=2l-s-t.
\end{array}
\end{eqnarray}

On the other hand, since the color degrees of the two vertices $u_0$
and $u_l$ are both at least $k$, and because of the assumption that
$P$ is one of the longest heterochromatic paths, we have that
$l+s\geq k$, $l+t\geq k$. This implies that $s\geq k-l$ and $t\geq
k-l$. Now we can get from Inequality \ref{ineq-3} that
$$k\leq 2l-s-t\leq 2l-2(k-l),$$
So, $4l\geq 3k$, and $l\geq \frac{3k}{4}$, and the proof is thus
complete. \qed

\section{Concluding remarks}

Finally, we give examples to show that our lower bounds given in
Theorems \ref{triangle-free-gallai-thm-1} and
\ref{triangle-free-gallai-thm-2} are best possible.

\begin{rem}
For any integer $k\geq 1$, there is a heterochromatic triangle free
complete graph $G_k$ with the color degree of every vertex $v$ in
$G_k$ is $k$, i.e., $d^c(v)=k$, such that any longest
heterochromatic $v$-path in $G_k$ is of length $k$.
\end{rem}

Let $G_k$ be an edge colored complete graph whose vertices are the
ordered $k$-tuples of 0's and 1's. An edge is in color $j$ ($1\leq
j\leq k$) if and only if the first $j-1$ coordinates of its two ends
are exactly the same and the j-th coordinates of its two ends are
different.

It is not hard to see that there exist no heterochromatic triangles
in $G_k$. Otherwise, suppose $uvw$ is a heterochromatic triangle,
where $u=(u_1,u_2,\ldots,u_k)$, $v=(v_1,v_2,\ldots,v_k)$ and
$w=(w_1,w_2,\ldots,w_k)$, the edge $uv$ is in color $x$, the edge
$vw$ is in color $y$, and the edge $uw$ is in color $z$. Without
loss of generality, we can assume that $1\leq x<y<z\leq k$. Since
the edge $uv$ is in color $x$, we can conclude that the first $x-1$
coordinates of the two vertices $u$ and $v$ are exactly the same,
and the $x$-th coordinates of $u$ and $v$ are different. Similarly,
we have that the first $z-1$ coordinates of the two vertices $u$ and
$w$ are exactly the same, and the $z$-th coordinates of $u$ and $w$
are different. Then the first $x-1$ coordinates of the two vertices
$v$ and $w$ are exactly the same, and the $x$-th coordinates of $v$
and $w$ are different. So the edge $vw$ is in color $x$, a
contradiction.

It is obvious that for every vertex $v$ in $G_k$, its color degree
is $k$. So we can easily conclude that the longest heterochromatic
$v$-path in $G_k$ is of length at least $k$ by Theorem
\ref{triangle-free-gallai-thm-1}. On the other hand, any longest
heterochromatic path in $G_k$ is not longer than $k$, since there
are only $k$ different colors used in this graph. Hence, the
conclusion in the remark is true.

\end{document}